\documentclass[11pt]{article}

\usepackage[utf8]{inputenc} 
\usepackage[T1]{fontenc}

\usepackage{amssymb, amsmath, amsthm, graphicx, enumerate, enumitem, color}
\usepackage{hyperref}
\usepackage{geometry}
\usepackage{algorithm}
\usepackage{algorithmic}
\usepackage[dvipsnames,svgnames,table]{xcolor}
\usepackage[normalem]{ulem}
\usepackage{booktabs}

\usepackage{etoolbox}

\usepackage{hyperref}
\usepackage[capitalise, compress, noabbrev]{cleveref}
\definecolor{linkblue}{named}{MidnightBlue}
\hypersetup{colorlinks=true, linkcolor=linkblue,  anchorcolor=linkblue,
	citecolor=linkblue, filecolor=linkblue, menucolor=linkblue,
	urlcolor=linkblue} 
\usepackage{mathtools}
\usepackage{thmtools, thm-restate}

\usepackage{footnote}
\makesavenoteenv{tabular}
\makesavenoteenv{table}

\usepackage{pgfplots}
\pgfplotsset{compat=1.18}

\theoremstyle{plain}
\newtheorem{thm}{Theorem}
\newtheorem*{thm*}{Theorem}
\let\leq\leqslant
\let\geq\geqslant

\let\epsilon\varepsilon

\let\setminus-

\DeclarePairedDelimiter\ceil{\lceil}{\rceil}
\DeclarePairedDelimiter\floor{\lfloor}{\rfloor}

\setenumerate{label=\textup{(\roman*)}, noitemsep, topsep=0pt,%-\parskip, 
labelindent=.2em, leftmargin=*, widest=iii,}
\setitemize{noitemsep, topsep=-\parskip, labelindent=.2em, leftmargin=*, widest=iii,}

\newcommand{\edge}[2]{#1 #2 }
\newcommand{\degree}[1]{\textrm{deg}({#1})}

\newcommand{\powerSum}[3]{\ensuremath{\sum_{i=#1}^{#2} (#3)^i}}
\newcommand{\rmg}{(\{r,m\};g)}

\newcommand{\numImprovedLowBounds}{123}

\title{Improved lower bounds on the maximum size of graphs with girth 5}

\providecommand{\keywords}[1]
{
  \small	
  \textbf{\textit{Keywords---}} #1
}

\date{}
\author{ 
Jan Goedgebeur \thanks{Department of Computer Science, KU Leuven Campus Kulak-Kortrijk, Kortrijk, Belgium. 
\protect\href{mailto:jan.goedgebeur@kuleuven.be}{\protect\nolinkurl{jan.goedgebeur@kuleuven.be}},
\protect\href{mailto:jorik.jooken@kuleuven.be}{\protect\nolinkurl{jorik.jooken@kuleuven.be}},
 and \protect\href{mailto:tibo.vandeneede@kuleuven.be}{\protect\nolinkurl{tibo.vandeneede@kuleuven.be}}} 
 \thanks{Department of Mathematics, Computer Science and Statistics, Ghent University, Ghent, Belgium.}
\and Jorik Jooken
 \footnotemark[1]
\and Gwenaël Joret \thanks{Computer Science Department, Université libre de Bruxelles, Brussels, Belgium. 
\protect\href{mailto:gwenael.joret@ulb.be}{\protect\nolinkurl{gwenael.joret@ulb.be}}} 
\and Tibo Van den Eede
 \footnotemark[1]
}

\begin{document}
\maketitle

\begin{abstract}
We present a new algorithm for improving lower bounds on
$ex(n;\{C_3,C_4\})$, the maximum size (number of edges) of an $n$-vertex graph of girth at least 5. 
The core of our algorithm is a variant of a hill-climbing heuristic introduced by Exoo, McKay, Myrvold and Nadon (2011) to find small cages. 
Our algorithm considers a range of values of $n$ in multiple passes.  
In each pass, the hill-climbing heuristic for a specific value of $n$ is initialized with a few graphs obtained by modifying near-extremal graphs previously found for neighboring values of $n$, allowing to `propagate' good patterns that were found. 
Focusing on the range $n\in \{74,75, \dots, 198\}$, which is currently beyond the scope of exact methods, 
our approach yields improvements on existing lower bounds for $ex(n;\{C_3,C_4\})$ for all $n$ in the range, except for two values of $n$ ($n=96,97$).  
\end{abstract}


\vspace*{1pt}

%\msc{05C07, 05C35, 05C85, 68R10, 90C35}

\section{Introduction}\label{sec:intro}
First, we introduce some basic terminology and notation. 
All graphs in this paper are simple and undirected. The set of vertices and set of edges of a graph $G$ are denoted by $V(G)$ and $E(G)$, respectively. 
The \textit{order} of a graph $G$ is its number of vertices $|V(G)|$, 
and its \textit{size} is its number of edges $|E(G)|$. 
An \textit{$s$-cycle} in a graph is a cycle with $s$ edges, and we denote by $C_s$ the graph which consists of an $s$-cycle.   
An edge between two vertices $v$ and $w$ is denoted by $\edge{v}{w}$. 
The degree of a vertex $v$ is denoted by $\degree{v}$. We call a vertex $v$ \textit{isolated} if $\degree{v} = 0$. A graph $G$ is \textit{$k$-regular} if each vertex of $G$ has degree $k$. The \textit{girth} of a graph $G$ is the length of the shortest cycle in $G$. If a graph has no cycle, then its girth is $\infty$.  
A graph $G$ is a \textit{subgraph} of a graph $H$ if $V(G) \subseteq V(H)$, $E(G) \subseteq E(H)$ and $E(H)$ only contains edges between vertices in $V(H)$. If $G$ is a subgraph of $H$, we also say that $H$ is a \textit{supergraph} of $G$. 

For a set $\mathcal{F}$ of graphs, $ex(n;\mathcal{F})$ denotes the maximum size in a graph of order $n$ that does not contain any member of $\mathcal{F}$ as a subgraph. The corresponding set of graphs of size $ex(n;\mathcal{F})$ is denoted by $EX(n;\mathcal{F})$ and are called \textit{extremal graphs}. 

The study of $ex(n;\mathcal{F})$ was initiated by Mantel in 1907~\cite{Mantel1907}, who proved that  $ex(n;\{C_3\}) = \floor{n^2/4}$, 
with the only extremal graph being the complete bipartite graph $K_{\floor{n/2},\ceil{n/2}}$. 
In 1941, Tur{\'a}n~\cite{Turan1941} subsequently showed that
\begin{equation*}
	ex(n;\{K_{r+1}\}) = \left( 1 - \frac{1}{r} + o(1) \right) \frac{n^2}{2}.
\end{equation*}
Erd\H{o}s and Stone~\cite{Erdos1946} generalized Tur{\'a}n's result in 1946, by showing that for every graph $H$ with chromatic number $r \geq 3$, we have
\begin{equation*}
	ex(n;\{H\}) = \left( \frac{r-2}{r-1} + o(1) \right) \binom{n}{2}.
\end{equation*}

Determining $ex(n;\mathcal{F})$ and $EX(n;\mathcal{F})$ for various sets of graphs $\mathcal{F}$ is a classical and much studied problem in extremal graph theory, see e.g.~\cite{AlonEtAl2003,ErdosSimonovits1966,Erdos1946,Furedi1983,Janzer2020,SudakovTomon2020}. 
In this paper, we focus on $ex(n;\{C_3,C_4\})$, which is thus the maximum size of a graph of order $n$ with girth at least 5. This problem was of particular interest to Erd\H{o}s. He conjectured in~\cite{Erdos1975} that \[ex(n;\{C_3,C_4\}) = \left(\frac{1}{2\sqrt{2}} + o(1) \right) n\sqrt{n}.\]
This conjecture is still open and has received significant attention over the years, see e.g.~\cite{AbajoEtAl2010new, AbajoEtAl2012girth, AbajoEtAl2010exact, AfzalyEtAl2023, FNV06, GarnickEtAl1993, GarnickEtAl1992, TangEtAl2009}. 
The current best asymptotic bounds are due to Garnick, Kwong and Lazebnik~\cite{GarnickEtAl1993}, who  proved that 
\begin{equation*}
    \frac{1}{2\sqrt{2}} \leq \limsup_{n\to\infty} \frac{ex(n;\{C_3,C_4\})}{n\sqrt{n}} \leq \frac{1}{2} .
\end{equation*}
We remark that, for $n\geq7$, it is known that every graph in $EX(n;\{C_3, C_4\})$  has girth exactly 5~\cite{GarnickEtAl1992}; thus,  $ex(n;\{C_3,C_4\})$ can be equivalently thought of as the maximum size of an $n$-vertex graph of girth $5$.  

While we focus on graphs with girth $5$ in this paper, we note that the corresponding extremal problem for graphs with larger girths (and related problems) have also been thoroughly studied in the literature, see e.g.~\cite{AbajoEtAl2010new,AbajoEtAl2012graphs, Aabajo2015, AbajoEtAl2010exact, Balbuena2008, PatternBoost, Marshall2011extremal,TangEtAl2009}. 

While Erd\H{o}s's conjecture on the asymptotic behavior of $ex(n;\{C_3,C_4\})$ is a central problem in the area, a second line of research that emerged over the years is the study of $ex(n;\{C_3,C_4\})$ for small values of $n$ using computer-assisted methods. The goal here is thus to obtain the best possible bounds on $ex(n;\{C_3,C_4\})$ for small $n$, typically up to around $200$ in the literature. This problem turned out to be a good benchmark for comparing different computational methods, and indeed a rich variety of methods have been applied to this problem over the years, such as greedy algorithms~\cite{AbajoEtAl2010exact,MarshallPhd2011}, hill-climbing with backtracking~\cite{GarnickEtAl1993}, a hybrid simulated annealing and genetic algorithm~\cite{TangEtAl2009}, tabu search~\cite{MehrabianEtAl2024}, and a reinforcement learning algorithm~\cite{MehrabianEtAl2024}. 
Nowadays, the exact value of $ex(n;\{C_3,C_4\})$ is known up to $n=53$~\cite{AfzalyEtAl2023}. For larger values of $n$, only lower bounds are known (see \cref{tab:summary_bounds}). 

In this paper, we introduce a new computational approach to find good lower bounds on $ex(n;\{C_3,C_4\})$. 
Our approach exploits a connection with the notorious \textit{Cage Problem}. Given two positive integers $k$ and $g$, the latter problem consists in determining the minimum order of a $k$-regular graph with girth exactly $g$. Let us call a \textit{$(k,g)$-graph} a $k$-regular graph with girth $g$, and a \textit{$(k,g)$-cage} one that has minimum order. A graph is a \textit{cage} if it is a $(k,g)$-cage for some integers $k,g$. A well-known lower bound on this minimum order is the Moore bound $M(k,g)$, defined as
% Our approach exploits a connection with the notorious \textit{Cage Problem}. Given two positive integers $k$ and $g$, the latter problem consists in determining the minimum order $n(k,g)$ of a $k$-regular graph with girth exactly $g$. Let us call \textit{$(k,g)$-graph} a $k$-regular graph with girth $g$, and \textit{$(k,g)$-cage} one that has minimum order (that is, order $n(k,g)$).  A well-known lower bound on $n(k,g)$ is the Moore bound $M(k,g)$, defined as
\begin{equation}\label{eq:Mbound}
	M(k,g) =
	\begin{cases}
		 1 + k \powerSum{0}{t-1}{k-1} & \text{for } g=2t+1 \\
		%\noalign{\vskip9pt}
		 2 \powerSum{0}{t-1}{k-1} & \text{for }g=2t . \\
	\end{cases}
\end{equation}
The $(k,g)$-cages that attain the Moore bound are called \textit{Moore graphs} and are also known as \textit{minimal cages}. 

The Cage Problem and the problem of determining $ex(n;\{C_3,C_4\})$ are different problems. Nevertheless, there is an interesting connection between the two problems, as shown by Abajo and Di{\'a}nez~\cite{AbajoEtAl2012graphs}.  

\begin{thm}[Theorem 4 in \cite{AbajoEtAl2012graphs}]\label{thm:mooreGraphsExtremal}
    Let $n$ and $g$ be positive integers with $g\geq 5$ and let $\mathcal{M}_{n,g}$ be the (possibly empty) set of Moore graphs of girth $g$ that have order $n$. If  $\mathcal{M}_{n,g}\neq \varnothing$, then 
    \begin{equation*}
		EX(n;\{C_3,\ldots, C_{g-1}\}) = \mathcal{M}_{n,g} .
	\end{equation*}
\end{thm}

In particular, the sizes of Moore graphs of girth $5$ determine $ex(n;\{C_3,C_4\})$ for the corresponding orders $n$. 
By the Hoffman-Singleton Theorem~\cite{HS60}, these graphs are respectively $C_5$ ($n=5$), the Petersen graph ($n=10$), the Hoffman-Singleton graph ($n=50$), and possibly some $57$-regular graph(s) of girth $5$ and order $n=3250$, the existence of which is a famous open problem. 

Another interesting observation in this direction is due to Backelin~\cite{Backelin2015}, who showed that for all $n\in \{40,45,47,48,49\}$, \emph{every} graph in $EX(n;\{C_3,C_{4}\})$ is in fact a subgraph of the Hoffman-Singleton graph (which is the unique graph in $EX(50;\{C_3,C_{4}\})$).

We see the above remarks as signs that cages of girth $5$ might be good candidates as seed graphs for local search heuristics aiming to find good lower bounds on $ex(n;\{C_3,C_4\})$ for the corresponding values of $n$. More generally, given a graph that gives a particularly good lower bound on $ex(n;\{C_3,C_4\})$ for some value of $n$, one might expect that some local modification of the graph could give good lower bounds  on $ex(n;\{C_3,C_4\})$ for nearby values of $n$ as well. This natural idea is heavily exploited in our algorithm: When trying to lower bound $ex(n;\{C_3,C_4\})$ for some $n$, we start with some well-chosen `good graphs' of order $n-1$ or $n+1$ that were previously found (i.e., graphs with girth at least 5 having many edges), modify them so that they have order $n$, and use the resulting graphs as seed graphs in a hill-climbing local search heuristic that tries to improve them as much as possible. The latter heuristic is a modification of an algorithm of Exoo, McKay, Myrvold and Nadon~\cite{ExooEtAl2011}, which they used for finding a $(4,7)$-cage.  

A detailed overview of our algorithm is given in \cref{sec:algo}. While the main ideas underlying our algorithm are relatively simple, there are several subtleties in its inner working. This is the result of a trial-and-error process in which we tried several variants of our algorithm before converging to the one presented in this paper, which is the version that gave the best results. 

Using our algorithm, we improved the best known lower bounds on $ex(n;\{C_3,C_4\})$ for all $n\in \{74,75, \dots, 198\}$, except for $n=96,97$ (where we matched the existing bounds). For some values of $n$, the improvements are particularly significant (in the double digits). 
See \cref{tab:summary_bounds} for a summary of the results.   
The corresponding graphs and the code of our implementation are available at \url{https://github.com/AGT-Kulak/searchMaxSizeMinGirth}~\cite{githubRepo}. 

We also mention that, as a by-product of the graphs found by our algorithm, we also obtained in passing four improved upper bounds on the minimum order of \textit{bi-regular cages} of girth $5$, a variant of cages where every vertex has degree in a prescribed set $\{r,m\}$; this is discussed \cref{sec:bi-regular cage}. 

The paper is organized as follows. In \cref{sec:algo}, we present our algorithm. In \cref{sec:results}, we discuss the results we obtained. Finally in \cref{sec:conc}, we give some concluding remarks and suggestions for future work.

\begin{table}[htbp]
	\centering
	\scalebox{0.88}{
		% \begin{tabular}{r r r || r r r || r r r || r r r }
                \begin{tabular}{c c c || c c c || c c c || c c c }
			\toprule
			% \multicolumn{1}{c}{$n$} & \multicolumn{2}{c}{$ex(n;\{C_3,C_4\})\geq$} & \multicolumn{1}{c}{$n$} & \multicolumn{2}{c}{$ex(n;\{C_3,C_4\})\geq$} & \multicolumn{1}{c}{$n$} & \multicolumn{2}{c}{$ex(n;\{C_3,C_4\})\geq$} & \multicolumn{1}{c}{$n$} & \multicolumn{2}{c}{$ex(n;\{C_3,C_4\})\geq$} \\
			% & Lit. & Alg. \ref{algo:compLowerBounds} & & Lit. & Alg. \ref{algo:compLowerBounds} & & Lit. & Alg. \ref{algo:compLowerBounds} & & Lit. & Alg. \ref{algo:compLowerBounds} \\ 
            \multicolumn{1}{c}{$n$} & \multicolumn{2}{c}{Lower bound} & \multicolumn{1}{c}{$n$} & \multicolumn{2}{c}{Lower bound} & \multicolumn{1}{c}{$n$} & \multicolumn{2}{c}{Lower bound} & \multicolumn{1}{c}{$n$} & \multicolumn{2}{c}{Lower bound} \\
			& Previous & New & & Previous & New & & Previous & New & & Previous & New  \\ 
			\midrule
			50 & 175 & 175 & 88 & 369 & \textbf{375} & 126 & 630 & \textbf{647} & 164 & 880 & \textbf{940} \\
			51 & 176 & 176 & 89 & 376 & \textbf{382} & 127 & 634 & \textbf{656} & 165 & 883 & \textbf{946} \\
			52 & 178 & 178 & 90 & 384 & \textbf{389} & 128 & 638 & \textbf{666} & 166 & 886 & \textbf{953} \\
			53 & 181 & 181 & 91 & 392 & \textbf{396} & 129 & 641 & \textbf{670} & 167 & 892 & \textbf{958} \\
			54 & 185 & 185 & 92 & 399 & \textbf{403} & 130 & 644 & \textbf{674} & 168 & 901 & \textbf{965} \\
			55 & 189 & 189 & 93 & 407 & \textbf{410} & 131 & 647 & \textbf{679} & 169 & 910 & \textbf{971} \\
			56 & 193 & 193 & 94 & 415 & \textbf{417} & 132 & 650 & \textbf{683} & 170 & 920 & \textbf{978} \\
			57 & 197 & 197 & 95 & 423 & \textbf{424} & 133 & 653 & \textbf{689} & 171 & 930 & \textbf{984} \\
			58 & 202 & 202 & 96 & 432 & 432 & 134 & 657 & \textbf{694} & 172 & 932 & \textbf{992} \\
			59 & 207 & 207 & 97 & 436 & 436 & 135 & 666 & \textbf{700} & 173 & 935 & \textbf{998} \\
			60 & 212 & 212 & 98 & 438 & \textbf{441} & 136 & 674 & \textbf{705} & 174 & 938 & \textbf{1005} \\
			61 & 216 & 216 & 99 & 440 & \textbf{446} & 137 & 683 & \textbf{711} & 175 & 941 & \textbf{1012} \\
			62 & 220 & 220 & 100 & 443 & \textbf{451} & 138 & 692 & \textbf{719} & 176 & 949 & \textbf{1020} \\
			63 & 224 & 224 & 101 & 445 & \textbf{457} & 139 & 700 & \textbf{727} & 177 & 958 & \textbf{1027} \\
			64 & 230 & 230 & 102 & 447 & \textbf{462} & 140 & 709 & \textbf{735} & 178 & 968 & \textbf{1035} \\
			65 & 235 & 235 & 103 & 452 & \textbf{468} & 141 & 717 & \textbf{743} & 179 & 977 & \textbf{1042} \\
			66 & 241 & 241 & 104 & 458 & \textbf{474} & 142 & 726 & \textbf{752} & 180 & 986 & \textbf{1050} \\
			67 & 246 & 246 & 105 & 464 & \textbf{482} & 143 & 735 & \textbf{760} & 181 & 995 & \textbf{1056} \\
			68 & 251 & 251 & 106 & 470 & \textbf{489} & 144 & 744 & \textbf{769} & 182 & 1004 & \textbf{1063} \\
			69 & 257 & 257 & 107 & 476 & \textbf{497} & 145 & 753 & \textbf{777} & 183 & 1013 & \textbf{1069} \\
			70 & 262 & 262 & 108 & 482 & \textbf{505} & 146 & 762 & \textbf{786} & 184 & 1022 & \textbf{1076} \\
			71 & 268 & 268 & 109 & 488 & \textbf{512} & 147 & 771 & \textbf{795} & 185 & 1032 & \textbf{1082} \\
			72 & 273 & 273 & 110 & 496 & \textbf{519} & 148 & 780 & \textbf{804} & 186 & 1042 & \textbf{1088} \\
			73 & 279 & 279 & 111 & 504 & \textbf{527} & 149 & 789 & \textbf{813} & 187 & 1052 & \textbf{1094} \\
			74 & 284 & \textbf{285} & 112 & 512 & \textbf{535} & 150 & 798 & \textbf{822} & 188 & 1062 & \textbf{1101} \\
			75 & 290 & \textbf{291} & 113 & 520 & \textbf{542} & 151 & 808 & \textbf{831} & 189 & 1072 & \textbf{1107} \\
			76 & 295 & \textbf{296} & 114 & 528 & \textbf{550} & 152 & 817 & \textbf{841} & 190 & 1082 & \textbf{1114} \\
			77 & 301 & \textbf{302} & 115 & 536 & \textbf{557} & 153 & 827 & \textbf{850} & 191 & 1092 & \textbf{1120} \\
			78 & 306 & \textbf{308} & 116 & 544 & \textbf{565} & 154 & 837 & \textbf{860} & 192 & 1102 & \textbf{1126} \\
			79 & 312 & \textbf{315} & 117 & 552 & \textbf{573} & 155 & 847 & \textbf{869} & 193 & 1112 & \textbf{1133} \\
			80 & 320 & \textbf{321} & 118 & 560 & \textbf{581} & 156 & 858 & \textbf{879} & 194 & 1122 & \textbf{1139} \\
			81 & 324 & \textbf{328} & 119 & 568 & \textbf{589} & 157 & 862 & \textbf{889} & 195 & 1132 & \textbf{1146} \\
			82 & 329 & \textbf{334} & 120 & 576 & \textbf{597} & 158 & 865 & \textbf{899} & 196 & 1142 & \textbf{1153} \\
			83 & 335 & \textbf{341} & 121 & 585 & \textbf{605} & 159 & 868 & \textbf{909} & 197 & 1152 & \textbf{1160} \\
			84 & 341 & \textbf{348} & 122 & 593 & \textbf{613} & 160 & 871 & \textbf{920} & 198 & 1163 & \textbf{1166} \\
			85 & 348 & \textbf{354} & 123 & 602 & \textbf{621} & 161 & 873 & \textbf{924} & & & \\
			86 & 355 & \textbf{361} & 124 & 611 & \textbf{629} & 162 & 875 & \textbf{929} & & & \\
			87 & 362 & \textbf{368} & 125 & 620 & \textbf{638} & 163 & 878 & \textbf{934} & & & \\
			\bottomrule
		\end{tabular}
	}
	\caption{Lower bounds on $ex(n;\{C_3,C_4\})$ for $50 \leq n \leq 198$: Best known bounds from the literature and bounds resulting from our algorithm; improved bounds are marked in bold. (Bibliographical references for previous lower bounds are given in \cref{sec:tables}.) }
	\label{tab:summary_bounds}
\end{table}

\section{Algorithm}\label{sec:algo}

The core of our algorithm is a local search heuristic that adds and deletes edges while preserving the girth constraint. This local search starts from good graphs already found for nearby values of $n$, that are first suitably modified.  
%and is initialized from small $(k,5)$–graphs and graphs of large size -- of which most are known to be extremal. 
\cref{subsec:localSearch} details the local search component that searches for graphs of large sizes, girth at least 5 and a specific order $n$.
While local search algorithms do not give any theoretical performance guarantees, they tend to yield very good results in practice at a scale where exact algorithms are no longer computationally feasible. In \cref{subsec:total_algo} we describe our full algorithm, which makes repeated uses of the local search as a subroutine over a given range of values for $n$. Each new local search begins with a graph that is a modification of a previously obtained graph of a nearby order, by adding or removing a vertex.

  %As start graphs we make use of small $(k,5)$-graphs and keep graphs of large size for later use. This will be outlined in \cref{subsec:total_algo}.

%% Some commands for variables in the local search algorithm
\newcommand{\localSearchAlgo}{\texttt{localSearch}}
\newcommand{\totalNumIters}{\texttt{totalNumIters}}
\newcommand{\currIter}{\texttt{currIter}}
\newcommand{\numItersTooRecent}{\texttt{numItersTooRecent}}
\newcommand{\maxk}{k_{\text{max}}}
\newcommand{\prob}{p}
\newcommand{\Gmax}{G_{\text{max}}}

\subsection{Local search}\label{subsec:localSearch}

In this subsection, we introduce a local search algorithm $\localSearchAlgo$ (\cref{algo:localSearch}) for finding graphs that have order~$n$, girth at least 5 and a large number of edges. The algorithm starts with an initial $n$-vertex graph $G$ of girth at least 5 and has the following  additional parameters:

\begin{itemize}
    \item $\totalNumIters$: the total number of iterations
    \item $\numItersTooRecent$: a threshold to prevent recently deleted edges from being removed again
     \item $\maxk$: the maximum number of edges to delete in one iteration     
    \item $\prob$: the probability of choosing an edge $\edge{u}{v}$ that maximizes $\degree{u} + \degree{v}$. 
\end{itemize}

The purpose of these parameters will be made clear later in the description of the algorithm (see \cref{algo:localSearch}). 
The algorithm consists of two nested loops: an outer loop, which repeats a number of local search trials, and an inner while loop that performs edge insertions until no more edges can be added without violating the girth constraint. 
As output we obtain a set of graphs $\mathcal{H}$ of order $n$ and girth at least 5, and sizes greater than or equal to that of $G$. All graphs in $\mathcal{H}$ have girth at least 5 since the starting graph $G\in \mathcal{H}$ has girth at least 5, and any other graph in $\mathcal{H}$ is obtained by the removal of edges from $G$ and addition of legal edges---which do not introduce 3-cycles nor 4-cycles.
An additional property of $\mathcal{H}$ is that no two graphs in $\mathcal{H}$ have the same sizes: 
a graph is added to $\mathcal{H}$ during the execution of the algorithm only if its size is greater than that of all graphs in $\mathcal{H}$. 
This might seem a peculiar choice, as opposed to e.g.\ outputting all the graphs of maximum size that were found (or of size at least $|E(G)|$), but we found that this strategy allows to keep $\mathcal{H}$ small while maintaining some diversity in $\mathcal{H}$, which leads to better results; this is discussed in \cref{subsec:total_algo}.

Our algorithm builds on the local search method by Exoo, McKay, Myrvold and Nadon~\cite{ExooEtAl2011}, 
which they used to find a $(4,7)$-cage. 
Given a seed graph, their algorithm first adds one by one \textit{legal edges} to the graph, edges that do not violate the girth constraint, each time selecting the legal edge with `highest priority', until no further legal edge can be added. At this point, some edges of the graph are deleted, based on a strategy that alternates between favoring older and newer edges, to enhance diversity in the search, and the whole process is started again. 

We made two main changes to their algorithm, which resulted in better bounds in our experiments. 
First, we changed the notion of priority for legal edges. 
The algorithm from~\cite{ExooEtAl2011} always chooses a legal edge $\edge{u}{v}$ with maximum degree sum $\degree{u}+\degree{v}$. 
In our algorithm, the choice of the next legal edge to be added is partly randomized:  
With probability $\prob$ we add a legal edge with maximum degree sum, and with probability $1-\prob$ we simply add a legal edge chosen uniformly at random. We tried different values between 0 and 1 for $\prob$, and a value around $0.5$ seems to give the best results. 
Second, in~\cite{ExooEtAl2011} edges are deleted either proportionally or inversely proportionally to their age, switching between these two methods after some number of iterations. 
In our experiments, we obtained better results simply by selecting some small number of edges to delete at random, and hence we chose this strategy.

\begin{algorithm}[ht!]
	\caption{\localSearchAlgo($G, \totalNumIters, \numItersTooRecent, \maxk, \prob$)}
	\label{algo:localSearch}
	\begin{algorithmic}[1]
		\STATE // $G$ is a graph of girth at least 5
		\STATE // $\totalNumIters$ is the number of iterations to perform
		\STATE // $\numItersTooRecent$ is a number of iterations such that the algorithm does not delete edges that were deleted in the last $\numItersTooRecent$ iterations
		\STATE // $\maxk$ is the maximum number of edges the algorithm can delete after an iteration
		\STATE // $\prob$ is the probability to choose the largest degree sum strategy when adding an edge		
		\STATE $\mathcal{H} \gets \{G\}$
        \STATE $c \gets |E(G)|$
		\WHILE{$\totalNumIters > 0$ and $c > 0$} 
		\WHILE{there exists a legal edge}
		\STATE Draw $u \sim \mathcal{U}(0,1)$
		\IF{$u < p$}
		\STATE Add a legal edge $\edge{u}{v}$ to $G$ with the largest degree sum $\degree{u}+\degree{v}$\\ 
        (in case of ties, choose a candidate pair uniformly  at random)
		\ELSE
		\STATE Add a legal edge to $G$, chosen uniformly at random
		\ENDIF
		\ENDWHILE
		 \IF{$|E(G)| > \max_{H \in \mathcal{H}} |E(H)|$}
		     \STATE $\mathcal{H} \gets \mathcal{H} \cup \{G\}$
		 \ENDIF	
            \STATE $c \gets $ number of edges of $G$ that have not been deleted in the last $\numItersTooRecent$ iterations
            \IF{$c \geq 1$}            
            \STATE Choose $k \in \{1,\ldots, \maxk\}$ uniformly at random
		\STATE Delete $\min\{k, c\}$ edges in $G$ that have not been deleted in the last $\numItersTooRecent$ iterations, chosen uniformly at random
            \ENDIF   
        \STATE $\totalNumIters \gets \totalNumIters - 1$
		\ENDWHILE		
		\RETURN $\mathcal{H}$
	\end{algorithmic}
\end{algorithm}

\subsection{Main algorithm}\label{subsec:total_algo}

\newcommand{\compLowerBoundsAlgo}{\texttt{computeLowerBounds}}
\newcommand{\downRunAlgo}{\texttt{downRun}}
\newcommand{\upRunAlgo}{\texttt{upRun}}

\newcommand{\bestGraphs}{\texttt{best}}
\newcommand{\emsGraphs}{\texttt{ems}}
\newcommand{\lsGraphs}{\texttt{ls}}

\newcommand{\nlow}{n_{\text{low}}}
\newcommand{\nup}{n_{\text{high}}}

In this subsection, we introduce our main algorithm, \cref{algo:compLowerBounds}, which computes lower bounds on $ex(n;\{C_3,C_4\})$ for all $n$ in some given interval $\{\nlow, \nlow+1, \dots,\nup\}$. 
The algorithm uses a data structure called $\bestGraphs$ to keep track of promising graphs found during the execution, that is, graphs with girth at least 5 and having a large number of edges. The graphs are stored per order, and we let  $\bestGraphs(n)$ denote the set of graphs in $\bestGraphs$ having order $n$. 
We initialize $\bestGraphs$ with the graphs from~\cite{AfzalyEtAl2023} (going up to order 64) and some small $(k,5)$-graphs (i.e., $(k,5)$-graphs with a small number of vertices), since they have a large size. Note that the use of small $(k,5)$-graphs is also motivated by \cref{thm:mooreGraphsExtremal} and the fact that the $(4,5)$- and $(6,5)$-cages,  which are not Moore graphs, are the only extremal graphs of their order. 
\cref{tab:initBest} gives a summary of the graphs that were used to initialize $\bestGraphs$. These graphs are also made available in \texttt{graph6} format in our GitHub repository~\cite{githubRepo}. 

\begin{table}[h]
	\centering
	\begin{tabular}{c c r}
		\toprule
		$n$ & \multicolumn{1}{c}{Comment} & \multicolumn{1}{c}{Reference} \\
		\midrule
		$50, \ldots, 53$ & All graphs from~\cite{AfzalyEtAl2023}, known to be extremal & \cite{AfzalyEtAl2023} \\
		$54, \ldots, 64$ & All graphs from~\cite{AfzalyEtAl2023}, not known to be extremal & \cite{AfzalyEtAl2023} \\
		$80$ & Smallest known $(8,5)$-graph & \cite{Royle2001} \\
		$96$ & Smallest known $(9,5)$-graph & \cite{Jorgensen2005} \\
		$124$ & Smallest known $(10,5)$-graph & \cite{ExooWeb} \\
		$126$ & $(10,5)$-graph & \cite{ExooWeb}  \\
		$154$ & Smallest known $(11,5)$-graph & \cite{ExooWeb}  \\
		$156$ & $(11,5)$-graph& \cite{ExooWeb}   \\
		$203$ & Smallest known $(12,5)$-graph & \cite{ExooWeb}  \\
		\bottomrule
	\end{tabular}
	\caption{The graphs used to initialize $\bestGraphs$ in \cref{algo:compLowerBounds}.
    } 
	\label{tab:initBest}
\end{table}

\cref{algo:compLowerBounds} iteratively applies two complementary routines: The ``up run'' (\cref{algo:upRun}) and the ``down run'' (\cref{algo:downRun}). Both of these routines update $\bestGraphs$ with the new graphs of large sizes found during their execution. 
For each $n \in \{\nlow, \ldots, \nup-1\}$ in increasing order, 
the up run starts with the top $\ell$ graphs of largest sizes in $\bestGraphs(n)$, adds an isolated vertex to each of them,  
and performs \cref{algo:localSearch} on the resulting graphs. 
The new graphs found by \cref{algo:localSearch} are then added to $\bestGraphs(n+1)$. 

The down run works in a similar fashion, considering each $n \in \{\nlow+1, \ldots, \nup\}$ in decreasing order this time, 
except the starting graphs are chosen as follows: 
First, we take the  top $\ell$ graphs in $\bestGraphs(n)$, then we consider all graphs obtained by deleting a vertex from these graphs in all possible ways, and finally we take the  top $\ell$ graphs in the resulting set of graphs.

Let us emphasize that $\upRunAlgo$ (\cref{algo:upRun}) and $\downRunAlgo$ (\cref{algo:downRun})  each perform two runs of the local search algorithm (\cref{algo:localSearch}) with different parameters. 
Since the local search algorithm is randomized, it usually makes sense to launch it more than once, to increase the probability of finding interesting graphs. 
In our implementation, we launch it only twice, to not increase the running time too much. We found that varying the parameters between the two runs give the best results as it promotes some variety in the resulting graphs. Lastly, as already mentioned in \cref{subsec:localSearch}, the local search only outputs graphs that have sizes greater than or equal to the size of the starting graph and that are all of different sizes. 
This might seem counter-intuitive at first sight---why not output, say, all found graphs of girth at least 5 with at least as many edges as the start graph?---but in our experiments we tried both approaches, and our approach with only saving graphs of strictly increasing size in $\mathcal{H}$ gave better final results.

We remark that in \cref{algo:upRun} and \cref{algo:downRun}, when adding graphs to the data structure $\bestGraphs$, we only add graphs that are not isomorphic to some graph already in $\bestGraphs$.  
This way $\bestGraphs$ is guaranteed to contain only non-isomorphic graphs, which is important for running time purposes, and for guaranteeing that when taking the top $\ell$ graphs of a given order, these are all non-isomorphic graphs. 
We use the \texttt{nauty} software library \cite{nauty} to efficiently filter isomorphic graphs.

\begin{algorithm}[ht!]
	\caption{\compLowerBoundsAlgo$(\nlow, \nup)$}
	\label{algo:compLowerBounds}
	\begin{algorithmic}[1]
		\STATE // $\nlow$ is the smallest order to compute lower bounds for
		\STATE // $\nup$ is the largest order to compute lower bounds for
            \STATE $\bestGraphs \gets$ graphs from \cref{tab:initBest}
		\WHILE{stopping criterion not met}
		\STATE $\bestGraphs \gets \upRunAlgo(\nlow, \nup, \bestGraphs,\ell=150)$ %// See \cref{algo:upRun}.
		\STATE $\bestGraphs \gets \downRunAlgo(\nlow, \nup, \bestGraphs,\ell=150)$ %// See \cref{algo:downRun}.
		\ENDWHILE
		\RETURN $\max_{G \in \bestGraphs(n)} |E(G)|$ for all orders $n \in \{\nlow, \ldots, \nup\}$
	\end{algorithmic}
\end{algorithm}

\begin{algorithm}[ht!]
	\caption{\upRunAlgo($\nlow, \nup, \bestGraphs, \ell$)}
	\label{algo:upRun}
	\begin{algorithmic}[1]
		\FOR{$n=\nlow$ up to $\nup-1$}
		\STATE Let $\mathcal{G}$ be the set of $\ell$ graphs with the most number of edges in $\bestGraphs(n)$
		\FOR{$G \in \mathcal{G}$}
		\STATE Let $G'$ be the graph obtained by adding an isolated vertex to $G$ 
            \STATE  Add all graphs in $\localSearchAlgo(G', 1000n, n, 3, 0.5)$ to $\bestGraphs(n+1)$ 
		\STATE  Add all graphs in $\localSearchAlgo(G', 1000n, \floor{n/3}, \max(3,\floor{n/10}), 0.5)$ to $\bestGraphs(n+1)$		
		\ENDFOR
		\ENDFOR
		\RETURN $\bestGraphs$
	\end{algorithmic}
\end{algorithm}

\begin{algorithm}[ht!]
	\caption{\downRunAlgo($\nlow, \nup, \bestGraphs, \ell$)}
	\label{algo:downRun}
	\begin{algorithmic}[1]
		\FOR{$n=\nup$ down to $\nlow+1$}
		\STATE $\mathcal{H} \gets \varnothing$
		\STATE Let $\mathcal{I}$ be the set of $\ell$ graphs with the most number of edges in $\bestGraphs(n)$
		\FOR{$G \in \mathcal{I}$}
		\FOR{$v \in V(G)$}
		\STATE $\mathcal{H} \gets  \mathcal{H} \cup \{G-v\}$
		\ENDFOR
		\ENDFOR
		\STATE Let $\mathcal{G}$ be the set of $\ell$ graphs with the most number of edges in $\mathcal{H}$
		\FOR{$G \in \mathcal{G}$}
		\STATE  Add all graphs in $\localSearchAlgo(G, 1000n, n, 3, 0.5)$ to $\bestGraphs(n-1)$ 
		\STATE  Add all graphs in $\localSearchAlgo(G, 1000n, \floor{n/3}, \max(3,\floor{n/10}), 0.5)$ to $\bestGraphs(n-1)$		
		\ENDFOR
		\ENDFOR
		\RETURN $\bestGraphs$
	\end{algorithmic}
\end{algorithm}

\section{Results}\label{sec:results}

In this section we present and discuss the results we obtained using  \cref{algo:compLowerBounds}.
We ran the algorithm for orders $n$ going from $\nlow=50$ up to $\nup=203$. 
We chose $\nlow=50$ because for this order it is known that the Hoffman-Singleton graph is extremal (c.f.\ the discussion around \cref{thm:mooreGraphsExtremal} in the introduction), and this is a few orders below $53$, the largest order for which $ex(n;\{C_3,C_4\})$ is known exactly. This way, the main algorithm is able to populate $\bestGraphs$ with some extra non-extremal graphs of order 53, which are then used in the searches for order 54 onwards. We chose $\nup=203$ because a small $(12,5)$-graph of that order is known, and this is also around the usual upper bounds on $n$ considered in the literature~\cite{ GarnickEtAl1993, MarshallPhd2011, MehrabianEtAl2024}. 

The parameter $\ell$ for both $\upRunAlgo$ and $\downRunAlgo$ and the parameter $\totalNumIters$ for $\localSearchAlgo$ were chosen to yield good results while keeping the CPU time within feasible limits. The remaining parameters of $\localSearchAlgo$ were determined experimentally by identifying the combinations that produced the best outcomes. 

We executed \cref{algo:compLowerBounds} for two full iterations of the \texttt{while} loop. We stopped after two iterations because the second iteration did not provide much improvement in comparison to the first one: For most values of $n$, there was no increase in the best lower bound. (For one value of $n$, there was an increase of $2$ though.) 
The two iterations required approximately $380$ CPU days in total. The top $\ell=150$ graphs of each order obtained after each up and down run can be found in \texttt{graph6} format at~\cite{githubRepo}. Additionally, a selection of the resulting graphs of largest size are available at the website ``House of Graphs'' \cite{HOG} by searching for the keyword ``maximum size of graphs with girth 5''.

\subsection{Improved bounds for $ex(n;\{C_3,C_4\})$}

The largest achieved sizes per order using \cref{algo:compLowerBounds} were already mentioned in the introduction, in \cref{tab:summary_bounds}. 
A more detailed overview of the results is given in  \cref{tab:results50to79,tab:results80to118,tab:results119to157,tab:results158to196,tab:results197to203} in \cref{sec:tables}, where the improvements made by each up and down run is described, and where bibliographical references are indicated for the best lower bounds reported in the literature\footnote{
We remark that we did not incorporate the lower bounds mentioned in~\cite{Bong2017} in the table because they are attributed to  an unpublished manuscript which does not seem to be available online, and neither the corresponding graphs nor the methods used to find them are described. 
Furthermore, we tried contacting the authors, without success. 
In any case, our results improve most of the bounds in~\cite{Bong2017}.}.  
Overall, our approach yields improved lower bounds for all $n\in \{74,75, \dots, 95\} \cup \{98, 99, \dots, 198\}$, thus for 
$\numImprovedLowBounds$ values of $n$, while for the remaining values we match the best known sizes from the literature. Recall the following result that was mentioned in the introduction due to Garnick, Kwong and Lazebnik~\cite{GarnickEtAl1993}: 
\begin{equation*}
    \frac{1}{2\sqrt{2}} \leq \limsup_{n\to\infty} \frac{ex(n;\{C_3,C_4\})}{n\sqrt{n}} \leq \frac{1}{2} .
\end{equation*}
In \cref{fig:plot} we show the new normalized lower bounds that we obtained.

\begin{figure}[t]
\centering
\begin{tikzpicture}
\begin{axis}[
    width=0.95\linewidth,
    height=6.2cm,
    xmin=50, xmax=198,
    ymin=0.34, ymax=0.51,
    xlabel={\(n\)},
    ylabel={\(\dfrac{\text{new lower bound}}{n^{3/2}}\)},
    tick align=outside,
    tick pos=left,
    grid=both,
    major grid style={gray!20},
    minor grid style={gray!10},
    minor tick num=1,
    xtick distance=10,
    ytick distance=0.02,
    legend style={
        draw=none,
        fill=none,
        font=\small,
        at={(0.5,1.08)},
        anchor=south,
    },
    legend columns=3,
    every axis plot/.style={line width=1.0pt},
    clip=false,
]
% main data
\addplot+[
    mark=*,
    mark size=1.3pt,
    color=blue!70!black,
] coordinates {
    (50,0.49497475)
(51,0.48323391)
(52,0.47469536)
(53,0.46909928)
(54,0.46620947)
(55,0.46335918)
(56,0.46054839)
(57,0.45777697)
(58,0.45730861)
(59,0.45676448)
(60,0.45615137)
(61,0.45337649)
(62,0.45064561)
(63,0.44795789)
(64,0.44921875)
(65,0.44843327)
(66,0.44947044)
(67,0.44856244)
(68,0.44762090)
(69,0.44839368)
(70,0.44735699)
(71,0.44796829)
(72,0.44685220)
(73,0.44732164)
(74,0.44771050)
(75,0.44802381)
(76,0.44675696)
(77,0.44696226)
(78,0.44710426)
(79,0.44861100)
(80,0.44861114)
(81,0.44993141)
(82,0.44980646)
(83,0.45095919)
(84,0.45202277)
(85,0.45172578)
(86,0.45264687)
(87,0.45349128)
(88,0.45426289)
(89,0.45496538)
(90,0.45560223)
(91,0.45617670)
(92,0.45669190)
(93,0.45715075)
(94,0.45755602)
(95,0.45791034)
(96,0.45927933)
(97,0.45638240)
(98,0.45456865)
(99,0.45277461)
(100,0.45100000)
(101,0.45022970)
(102,0.44847859)
(103,0.44770301)
(104,0.44691850)
(105,0.44798460)
(106,0.44807433)
(107,0.44903555)
(108,0.44994118)
(109,0.44991473)
(110,0.44986099)
(111,0.45063627)
(112,0.45136382)
(113,0.45121302)
(114,0.45186158)
(115,0.45165668)
(116,0.45223218)
(117,0.45276811)
(118,0.45326598)
(119,0.45372724)
(120,0.45415329)
(121,0.45454545)
(122,0.45490502)
(123,0.45523321)
(124,0.45553119)
(125,0.45651564)
(126,0.45745509)
(127,0.45835108)
(128,0.45989562)
(129,0.45728846)
(130,0.45472085)
(131,0.45219156)
(132,0.45189088)
(133,0.45157559)
(134,0.45038189)
(135,0.44961004)
(136,0.44918563)
(137,0.44879913)
(138,0.44844941)
(139,0.44813538)
(140,0.44832825)
(141,0.44803544)
(142,0.44777608)
(143,0.44754878)
(144,0.44735214)
(145,0.44718479)
(146,0.44704534)
(147,0.44693236)
(148,0.44684451)
(149,0.44678043)
(150,0.44673878)
(151,0.44671819)
(152,0.44720186)
(153,0.44644553)
(154,0.44630105)
(155,0.44599486)
(156,0.44555405)
(157,0.44499336)
(158,0.44433564)
(159,0.44360229)
(160,0.44341497)
(161,0.44319018)
(162,0.44293676)
(163,0.44266368)
(164,0.43925207)
(165,0.43744595)
(166,0.43568376)
(167,0.43447538)
(168,0.43327549)
(169,0.43185044)
(170,0.43082878)
(171,0.42983496)
(172,0.42915914)
(173,0.42849633)
(174,0.42784552)
(175,0.42720576)
(176,0.42657618)
(177,0.42595604)
(178,0.42534456)
(179,0.42474110)
(180,0.42414500)
(181,0.42355571)
(182,0.42297265)
(183,0.42239531)
(184,0.42182324)
(185,0.42125615)
(186,0.42069367)
(187,0.42013562)
(188,0.42781392)
(189,0.42604445)
(190,0.42535818)
(191,0.42429506)
(192,0.42323985)
(193,0.42256541)
(194,0.42152286)
(195,0.42085520)
(196,0.42018950)
(197,0.41952578)
(198,0.41850511)
};
\addlegendentry{new lower bound}

% benchmark 1/2^{3/2}
\addplot+[
    very thick,
    red,
    solid,
    mark=none,
] coordinates {
    (50,0.35355339)
    (198,0.35355339)
};
\addlegendentry{$1/2^{3/2}$}

% benchmark 1/2
\addplot+[
    very thick,
    black,
    solid,
    mark=none,
] coordinates {
    (50,0.5)
    (198,0.5)
};
\addlegendentry{$1/2$}

\end{axis}
\end{tikzpicture}
\caption{Normalized lower bounds from Table~\ref{tab:summary_bounds}.}
\label{fig:plot}
\end{figure}

% Furthermore, all graphs that we obtained of largest size for its order have diameter 3, except for the Hoffman-Singleton graph of order 50 which has diameter 2.

Here are some comments and observations about the obtained lower bounds and the corresponding graphs.
First, we do not obtain improvements for $n \in \{54, \ldots, 73\} \cup \{96, 97\} \cup\{199, \ldots, 203\}$. We believe that the existing lower bounds for these orders are already (quasi-)optimal. The smaller orders $n \in \{54, \ldots, 73\}$, in particular, have been studied more extensively, for instance by Afzaly and McKay~\cite{AfzalyEtAl2023} up to $n=64$, who also determined $ex(n;\{C_3,C_4\})$ exactly up to $n=53$. For order $96$ and $203$ we note that the lower bound corresponds to the size of the smallest known $(9,5)$- and $(12,5)$-graph respectively. Marshall~\cite{MarshallPhd2011} also exploited these graphs for their own and nearby orders, albeit with a different, greedy algorithm.  
The $(9,5)$- and $(12,5)$-graphs might be extremal, just as the $(k,5)$-cages for $k \in \{2,3,4,6,7\}$ are. 

On the other hand, and perhaps surprisingly, we also found that the smallest known $(k,5)$-graphs are not always the best graphs for their respective orders:

\begin{itemize}
\item $ex(80;\{C_3,C_4\})\geq 321$, while the smallest known $(8,5)$-graph has order 80 and size 320;
\item $ex(124;\{C_3,C_4\})\geq 629$, while the smallest known $(10,5)$-graph has order 124 and size 620;
\item $ex(154;\{C_3,C_4\})\geq 860$, while the smallest known $(11,5)$-graph has order 154 and size 847.
\end{itemize}

$ $

Finally, we remark that for certain values of $n$, our lower bound improvements are remarkably large. 
For instance, for both $n=175$ and $n=176$, we increased the existing lower bound by $71$. 
We also observed that the graphs found for the range of values of $n$ we considered tend to have many vertex orbits and we were unable to detect any clear patterns in them that would suggest good constructions for larger values of $n$.

\subsection{Improved bounds on the bi-regular cage problem}
\label{sec:bi-regular cage}

In the introduction we mentioned a connection between the classical Cage Problem and the problem of determining $ex(n;\{C_3, \ldots, C_{g-1}\})$. 
Given positive integers $r,m$ and $g$ with $r<m$ and $g\geq 3$, a variant of the Cage Problem is to search for $\rmg$-graphs of minimum order, where an \textit{$\rmg$-graph} is a graph of girth $g$ with degree set $\{r,m\}$. 
An \textit{$\rmg$-cage}, also known as a \textit{bi-regular cage}, is an $\rmg$-graph of minimum order. 
The order of an $\rmg$-cage is denoted by $n\rmg$.

As it happens, some of the graphs found by our algorithm also improve some of the best known upper bounds on $n\rmg$. For an overview of the best known bounds on $n\rmg$ for $r\leq5$, we refer the reader to~\cite{GoedgebeurEtAl2024}.~\cref{tab:biregCages} summarizes the improvements on upper bounds for $n\rmg$ that we obtained. The corresponding graphs are also available at~\cite{githubRepo}.

\begin{table}[h]
	\centering
	\begin{tabular}{r r r r r}
		\toprule
		$r$ & $m$ & \multicolumn{2}{c}{$n(\{r,m\};5)\leq$} & Source literature \\
		\cmidrule(lr){3-4}
		&  & Literature &  \cref{algo:compLowerBounds}& \\
		\midrule
		8 & 9 & 96  & \textbf{88} &  Prop. 2 (a)~\cite{Dn-cages}\\
		9 & 12 & 193  & \textbf{180} & Th. 2~\cite{Yang2003} + Th. 2.1~\cite{monotoDn-cages}\footnote{Th.~2~\cite{Yang2003} shows that $n(\{9,12\};6)\leq194$ and Th.~2.1~\cite{monotoDn-cages} shows that $n\rmg < n(\{r,m\};g+1)$. Hence, $n(\{9,12\};5)\leq193$.} \\
		10 & 11 & 154  & \textbf{128} & Prop. 2 (a)~\cite{Dn-cages}\\
		11 & 12 & 203  & \textbf{155} & Prop. 2 (a)~\cite{Dn-cages}\\
		\bottomrule
	\end{tabular}
	\caption{Improved upper bounds on $n(\{r,m\};5)$ by \cref{algo:compLowerBounds} (marked in bold) in comparison to the literature.
    }
	\label{tab:biregCages}
\end{table}

\section{Conclusion}\label{sec:conc}

In this work, we introduced an algorithm that searches for graphs with given order $n$, girth at least 5, and having as many edges as possible, 
where $n$ belongs to some fixed range of values. For a given order $n$, the algorithm uses a carefully crafted local search heuristic, which is a modification of an algorithm due to Exoo, McKay, Myrvold and Nadon~\cite{ExooEtAl2011} for finding small cages. 
Each local search is started from some promising graphs obtained from the best graphs found so far for nearby orders ($n-1$ and $n+1$), and also (crucially) from some known $(k,5)$-graphs of small order.  Applying this method on the range $n\in \{50, 51, \dots, 203\}$, 
we improved the existing lower bounds on $ex(n;\{C_3,C_4\})$ for all $n\in \{74,75, \dots, 95\} \cup \{98, 99, \dots, 198\}$, thus for 
$\numImprovedLowBounds$ values of $n$, while for the remaining values we match the best known sizes from the literature. Some of the improvements are relatively large: For example, for $n=175,176$, we increased the lower bound by 71 edges.

As a by-product of our methods, we also obtained four improved upper bounds on $n(\{r,m\};5)$, which is the minimum order of a graph of girth 5 and degree set $\{r,m\}$ with $r<m$. 

For future works, a natural direction would be to experiment further with the general approach of \cref{algo:compLowerBounds} but replace our local search heuristic, which is used as a black box, by another heuristic. It could very well be the case that a completely different heuristic gives even better results. We remark that some improvements might also be achievable simply by spending more CPU time on the execution of the algorithm, though we expect it would only increase a few lower bounds by some $+1$s or $+2$s, thus we do not see this as a promising direction. We believe that new ideas are necessary for further significant improvements. We also believe that the lower bounds are now very close to the exact value of $ex(n;\{C_3,C_4\})$ for, say, $n\leq 100$, while there might still be room for improvements for larger $n$. The belief is based on the observation that extremal and near-extremal graphs of a given order are often closely related to extremal or near-extremal graphs of nearby orders, as illustrated by \cite{AfzalyEtAl2023,Backelin2015}. In particular, the graphs known to be extremal up to order $53$, together with the graphs of slightly larger order that are already close to the best known upper bounds, suggest that the graphs produced by our algorithm---obtained by increasing the order of such seed graphs and then making local modifications---are likely to remain near-extremal.

A second direction for future work would be to consider the same problem but for graphs of larger girths. Indeed, it is straightforward to adapt our algorithm to find lower bounds on $ex(n;\{C_3,\ldots,C_{g-1}\})$ with $g\geq 6$. This is something we have not explored at all. 

Finally, it would be very interesting to identify patterns in some of the graphs found by our algorithm, which would suggest good constructions for larger $n$'s. However, this is probably difficult because most of these graphs have no automorphisms except for the trivial one.

\section*{Acknowledgements}
Jan Goedgebeur and Tibo Van den Eede are supported by a grant of the Research Foundation Flanders (FWO) with grant number G0AGX24N and by Internal Funds of KU Leuven. Jorik Jooken is supported by an FWO grant with grant number 1222524N. 
Gwenaël Joret is supported by the Fonds National de la Recherche Scientifique (F.R.S.--FNRS). 
Gwenaël Joret thanks Adam Zsolt Wagner for introducing him to the problem studied in this paper and for enlightening discussions.  

\bibliographystyle{abbrv}
\typeout{}
\bibliography{references}

\newpage

\appendix

\section{Lower bounds for $ex(n;\{C_3,C_4\})$ }\label{sec:tables}

\begin{table}[h]
\centering
\begin{tabular}{r r r |  r r r r }
\toprule
\multicolumn{1}{c}{$n$} & \multicolumn{6}{c}{$ex(n;\{C_3,C_4\})\geq$} \\
\cmidrule(lr){2-7}
 & & & \multicolumn{4}{c}{\cref{algo:compLowerBounds}} \\ 
 \cmidrule(lr){4-7}
 & Literature & Source & Up 1 & Down 1 & Up 2 & Down 2 \\ 
\midrule
50 & 175 & \cite{AfzalyEtAl2023} & - & 175 & - & 175 \\
51 & 176 & \cite{AfzalyEtAl2023} & 176 & 176 & 176 & 176 \\
52 & 178 & \cite{AfzalyEtAl2023} & 178 & 178 & 178 & 178 \\
53 & 181 & \cite{AfzalyEtAl2023} & 181 & 181 & 181 & 181 \\
54 & 185 & \cite{GarnickEtAl1993} & 185 & 185 & 185 & 185 \\
55 & 189 & \cite{AfzalyEtAl2023} & 189 & 189 & 189 & 189 \\
56 & 193 & \cite{AfzalyEtAl2023} & 193 & 193 & 193 & 193 \\
57 & 197 & \cite{AfzalyEtAl2023} & 197 & 197 & 197 & 197 \\
58 & 202 & \cite{AfzalyEtAl2023} & 202 & 202 & 202 & 202 \\
59 & 207 & \cite{AfzalyEtAl2023} & 207 & 207 & 207 & 207 \\
60 & 212 & \cite{AfzalyEtAl2023} & 212 & 212 & 212 & 212 \\
61 & 216 & \cite{AfzalyEtAl2023} & 216 & 216 & 216 & 216 \\
62 & 220 & \cite{AfzalyEtAl2023} & 220 & 220 & 220 & 220 \\
63 & 224 & \cite{AfzalyEtAl2023} & 224 & 224 & 224 & 224 \\
64 & 230 & \cite{AfzalyEtAl2023} & 230 & 230 & 230 & 230 \\
65 & 235 & \cite{MehrabianEtAl2024} & 235 & 235 & 235 & 235 \\
66 & 241 & \cite{MehrabianEtAl2024} & 240 & 241 & 241 & 241 \\
67 & 246 & \cite{MehrabianEtAl2024} & 246 & 246 & 246 & 246 \\
68 & 251 & \cite{MehrabianEtAl2024} & 251 & 251 & 251 & 251 \\
69 & 257 & \cite{MehrabianEtAl2024} & 257 & 257 & 257 & 257 \\
70 & 262 & \cite{MehrabianEtAl2024} & 262 & 262 & 262 & 262 \\
71 & 268 & \cite{MehrabianEtAl2024} & 268 & 268 & 268 & 268 \\
72 & 273 & \cite{MehrabianEtAl2024} & 273 & 273 & 273 & 273 \\
73 & 279 & \cite{MehrabianEtAl2024} & 279 & 279 & 279 & 279 \\
74 & 284 & \cite{MehrabianEtAl2024} & \textbf{285} & \textbf{285} & \textbf{285} & \textbf{285} \\
75 & 290 & \cite{MehrabianEtAl2024} & \textbf{291} & \textbf{291} & \textbf{291} & \textbf{291} \\
76 & 295 & \cite{MehrabianEtAl2024} & \textbf{296} & \textbf{296} & \textbf{296} & \textbf{296} \\
77 & 301 & \cite{MehrabianEtAl2024} & \textbf{302} & \textbf{302} & \textbf{302} & \textbf{302} \\
78 & 306 & \cite{MehrabianEtAl2024} & \textbf{307} & \textbf{308} & \textbf{308} & \textbf{308} \\
79 & 312 & \cite{MarshallPhd2011} & \textbf{313} & \textbf{315} & \textbf{315} & \textbf{315} \\
\bottomrule
\end{tabular}
\caption{Lower bounds of $ex(n;\{C_3,C_4\})$ for $50 \leq n \leq 79$ from the literature and \cref{algo:compLowerBounds}. Bounds which improve upon the ones from the literature are marked in bold.}
\label{tab:results50to79}
\end{table}

\begin{table}[h]
\centering
\begin{tabular}{r r r |  r r r r }
\toprule
\multicolumn{1}{c}{$n$} & \multicolumn{6}{c}{$ex(n;\{C_3,C_4\})\geq$} \\
\cmidrule(lr){2-7}
 & & & \multicolumn{4}{c}{\cref{algo:compLowerBounds}} \\ 
 \cmidrule(lr){4-7}
 & Literature & Source & Up 1 & Down 1 & Up 2 & Down 2 \\ 
\midrule
80 & 320 & \cite{MarshallPhd2011} & 320 & \textbf{321} & \textbf{321} & \textbf{321} \\
81 & 324 & \cite{MehrabianEtAl2024} & 324 & \textbf{328} & \textbf{328} & \textbf{328} \\
82 & 329 & \cite{MehrabianEtAl2024} & \textbf{330} & \textbf{334} & \textbf{334} & \textbf{334} \\
83 & 335 & \cite{MehrabianEtAl2024} & 335 & \textbf{341} & \textbf{341} & \textbf{341} \\
84 & 341 & \cite{MarshallPhd2011} & 341 & \textbf{348} & \textbf{348} & \textbf{348} \\
85 & 348 & \cite{MarshallPhd2011} & 346 & \textbf{354} & \textbf{354} & \textbf{354} \\
86 & 355 & \cite{MarshallPhd2011} & 352 & \textbf{361} & \textbf{361} & \textbf{361} \\
87 & 362 & \cite{MarshallPhd2011} & 358 & \textbf{368} & \textbf{368} & \textbf{368} \\
88 & 369 & \cite{MarshallPhd2011} & 363 & \textbf{375} & \textbf{375} & \textbf{375} \\
89 & 376 & \cite{MarshallPhd2011} & 369 & \textbf{382} & \textbf{382} & \textbf{382} \\
90 & 384 & \cite{MarshallPhd2011} & 375 & \textbf{389} & \textbf{389} & \textbf{389} \\
91 & 392 & \cite{MarshallPhd2011} & 382 & \textbf{396} & \textbf{396} & \textbf{396} \\
92 & 399 & \cite{MarshallPhd2011} & 387 & \textbf{403} & \textbf{403} & \textbf{403} \\
93 & 407 & \cite{MarshallPhd2011} & 394 & \textbf{410} & \textbf{410} & \textbf{410} \\
94 & 415 & \cite{MarshallPhd2011} & 400 & \textbf{417} & \textbf{417} & \textbf{417} \\
95 & 423 & \cite{MarshallPhd2011} & 405 & \textbf{424} & \textbf{424} & \textbf{424} \\
96 & 432 & \cite{MarshallPhd2011} & 432 & 432 & 432 & 432 \\
97 & 436 & \cite{MarshallPhd2011} & 436 & 436 & 436 & 436 \\
98 & 438 & \cite{MarshallPhd2011} & \textbf{441} & \textbf{441} & \textbf{441} & \textbf{441} \\
99 & 440 & \cite{MarshallPhd2011} & \textbf{446} & \textbf{446} & \textbf{446} & \textbf{446} \\
100 & 443 & \cite{MarshallPhd2011} & \textbf{451} & \textbf{451} & \textbf{451} & \textbf{451} \\
101 & 445 & \cite{MarshallPhd2011} & \textbf{457} & \textbf{457} & \textbf{457} & \textbf{457} \\
102 & 447 & \cite{MarshallPhd2011} & \textbf{462} & \textbf{462} & \textbf{462} & \textbf{462} \\
103 & 452 & \cite{MehrabianEtAl2024} & \textbf{468} & \textbf{468} & \textbf{468} & \textbf{468} \\
104 & 458 & \cite{MehrabianEtAl2024} & \textbf{474} & \textbf{474} & \textbf{474} & \textbf{474} \\
105 & 464 & \cite{MehrabianEtAl2024} & \textbf{480} & \textbf{482} & \textbf{482} & \textbf{482} \\
106 & 470 & \cite{MehrabianEtAl2024} & \textbf{487} & \textbf{489} & \textbf{489} & \textbf{489} \\
107 & 476 & \cite{MehrabianEtAl2024} & \textbf{492} & \textbf{497} & \textbf{497} & \textbf{497} \\
108 & 482 & \cite{MehrabianEtAl2024} & \textbf{498} & \textbf{504} & \textbf{504} & \textbf{505} \\
109 & 488 & \cite{MarshallPhd2011} & \textbf{505} & \textbf{511} & \textbf{512} & \textbf{512} \\
110 & 496 & \cite{MarshallPhd2011} & \textbf{511} & \textbf{519} & \textbf{519} & \textbf{519} \\
111 & 504 & \cite{MarshallPhd2011} & \textbf{517} & \textbf{526} & \textbf{527} & \textbf{527} \\
112 & 512 & \cite{MarshallPhd2011} & \textbf{524} & \textbf{533} & \textbf{535} & \textbf{535} \\
113 & 520 & \cite{MarshallPhd2011} & \textbf{531} & \textbf{541} & \textbf{542} & \textbf{542} \\
114 & 528 & \cite{MarshallPhd2011} & \textbf{537} & \textbf{549} & \textbf{549} & \textbf{550} \\
115 & 536 & \cite{MarshallPhd2011} & \textbf{544} & \textbf{557} & \textbf{557} & \textbf{557} \\
116 & 544 & \cite{MarshallPhd2011} & \textbf{551} & \textbf{564} & \textbf{564} & \textbf{565} \\
117 & 552 & \cite{MarshallPhd2011} & \textbf{557} & \textbf{572} & \textbf{572} & \textbf{573} \\
118 & 560 & \cite{MarshallPhd2011} & \textbf{564} & \textbf{580} & \textbf{580} & \textbf{581} \\
\bottomrule
\end{tabular}
\caption{Lower bounds of $ex(n;\{C_3,C_4\})$ for $80 \leq n \leq 118$ from the literature and \cref{algo:compLowerBounds}. Bounds which improve upon the ones from the literature are marked in bold.}
\label{tab:results80to118}
\end{table}

\begin{table}[h]
\centering
\begin{tabular}{r r r |  r r r r }
\toprule
\multicolumn{1}{c}{$n$} & \multicolumn{6}{c}{$ex(n;\{C_3,C_4\})\geq$} \\
\cmidrule(lr){2-7}
 & & & \multicolumn{4}{c}{\cref{algo:compLowerBounds}} \\ 
 \cmidrule(lr){4-7}
 & Literature & Source & Up 1 & Down 1 & Up 2 & Down 2 \\ 
\midrule
119 & 568 & \cite{MarshallPhd2011} & \textbf{570} & \textbf{588} & \textbf{588} & \textbf{589} \\
120 & 576 & \cite{MarshallPhd2011} & \textbf{577} & \textbf{597} & \textbf{597} & \textbf{597} \\
121 & 585 & \cite{MarshallPhd2011} & 582 & \textbf{604} & \textbf{605} & \textbf{605} \\
122 & 593 & \cite{MarshallPhd2011} & 588 & \textbf{613} & \textbf{613} & \textbf{613} \\
123 & 602 & \cite{MarshallPhd2011} & 594 & \textbf{621} & \textbf{621} & \textbf{621} \\
124 & 611 & \cite{MarshallPhd2011} & \textbf{620} & \textbf{629} & \textbf{629} & \textbf{629} \\
125 & 620 & \cite{MarshallPhd2011} & \textbf{637} & \textbf{638} & \textbf{638} & \textbf{638} \\
126 & 630 & \cite{MarshallPhd2011} & \textbf{647} & \textbf{647} & \textbf{647} & \textbf{647} \\
127 & 634 & \cite{MarshallPhd2011} & \textbf{656} & \textbf{656} & \textbf{656} & \textbf{656} \\
128 & 638 & \cite{MarshallPhd2011} & \textbf{666} & \textbf{666} & \textbf{666} & \textbf{666} \\
129 & 641 & \cite{MarshallPhd2011} & \textbf{669} & \textbf{670} & \textbf{670} & \textbf{670} \\
130 & 644 & \cite{MarshallPhd2011} & \textbf{673} & \textbf{673} & \textbf{674} & \textbf{674} \\
131 & 647 & \cite{MarshallPhd2011} & \textbf{678} & \textbf{679} & \textbf{679} & \textbf{679} \\
132 & 650 & \cite{MarshallPhd2011} & \textbf{683} & \textbf{683} & \textbf{683} & \textbf{683} \\
133 & 653 & \cite{MarshallPhd2011} & \textbf{688} & \textbf{689} & \textbf{689} & \textbf{689} \\
134 & 657 & \cite{MarshallPhd2011} & \textbf{694} & \textbf{694} & \textbf{694} & \textbf{694} \\
135 & 666 & \cite{MarshallPhd2011} & \textbf{699} & \textbf{699} & \textbf{699} & \textbf{700} \\
136 & 674 & \cite{MarshallPhd2011} & \textbf{705} & \textbf{705} & \textbf{705} & \textbf{705} \\
137 & 683 & \cite{MarshallPhd2011} & \textbf{710} & \textbf{710} & \textbf{711} & \textbf{711} \\
138 & 692 & \cite{MarshallPhd2011} & \textbf{716} & \textbf{719} & \textbf{719} & \textbf{719} \\
139 & 700 & \cite{MarshallPhd2011} & \textbf{722} & \textbf{727} & \textbf{727} & \textbf{727} \\
140 & 709 & \cite{MarshallPhd2011} & \textbf{728} & \textbf{735} & \textbf{735} & \textbf{735} \\
141 & 717 & \cite{MarshallPhd2011} & \textbf{734} & \textbf{743} & \textbf{743} & \textbf{743} \\
142 & 726 & \cite{MarshallPhd2011} & \textbf{740} & \textbf{752} & \textbf{752} & \textbf{752} \\
143 & 735 & \cite{MarshallPhd2011} & \textbf{746} & \textbf{760} & \textbf{760} & \textbf{760} \\
144 & 744 & \cite{MarshallPhd2011} & \textbf{752} & \textbf{769} & \textbf{769} & \textbf{769} \\
145 & 753 & \cite{MarshallPhd2011} & \textbf{758} & \textbf{777} & \textbf{777} & \textbf{777} \\
146 & 762 & \cite{MarshallPhd2011} & \textbf{764} & \textbf{786} & \textbf{786} & \textbf{786} \\
147 & 771 & \cite{MarshallPhd2011} & 771 & \textbf{795} & \textbf{795} & \textbf{795} \\
148 & 780 & \cite{MarshallPhd2011} & 777 & \textbf{804} & \textbf{804} & \textbf{804} \\
149 & 789 & \cite{MarshallPhd2011} & 784 & \textbf{813} & \textbf{813} & \textbf{813} \\
150 & 798 & \cite{MarshallPhd2011} & 790 & \textbf{822} & \textbf{822} & \textbf{822} \\
151 & 808 & \cite{MarshallPhd2011} & 797 & \textbf{831} & \textbf{831} & \textbf{831} \\
152 & 817 & \cite{MarshallPhd2011} & 803 & \textbf{841} & \textbf{841} & \textbf{841} \\
153 & 827 & \cite{MarshallPhd2011} & 810 & \textbf{850} & \textbf{850} & \textbf{850} \\
154 & 837 & \cite{MarshallPhd2011} & \textbf{847} & \textbf{860} & \textbf{860} & \textbf{860} \\
155 & 847 & \cite{MarshallPhd2011} & \textbf{866} & \textbf{869} & \textbf{869} & \textbf{869} \\
156 & 858 & \cite{MarshallPhd2011} & \textbf{877} & \textbf{879} & \textbf{879} & \textbf{879} \\
157 & 862 & \cite{MarshallPhd2011} & \textbf{888} & \textbf{889} & \textbf{889} & \textbf{889} \\
\bottomrule
\end{tabular}
\caption{Lower bounds of $ex(n;\{C_3,C_4\})$ for $119 \leq n \leq 157$ from the literature and \cref{algo:compLowerBounds}. Bounds which improve upon the ones from the literature are marked in bold.}
\label{tab:results119to157}
\end{table}

\begin{table}[h]
\centering
\begin{tabular}{r r r |  r r r r }
\toprule
\multicolumn{1}{c}{$n$} & \multicolumn{6}{c}{$ex(n;\{C_3,C_4\})\geq$} \\
\cmidrule(lr){2-7}
 & & & \multicolumn{4}{c}{\cref{algo:compLowerBounds}} \\ 
 \cmidrule(lr){4-7}
 & Literature & Source & Up 1 & Down 1 & Up 2 & Down 2 \\ 
\midrule
158 & 865 & \cite{MarshallPhd2011} & \textbf{898} & \textbf{899} & \textbf{899} & \textbf{899} \\
159 & 868 & \cite{MarshallPhd2011} & \textbf{909} & \textbf{909} & \textbf{909} & \textbf{909} \\
160 & 871 & \cite{MarshallPhd2011} & \textbf{920} & \textbf{920} & \textbf{920} & \textbf{920} \\
161 & 873 & \cite{MarshallPhd2011} & \textbf{924} & \textbf{924} & \textbf{924} & \textbf{924} \\
162 & 875 & \cite{MarshallPhd2011} & \textbf{929} & \textbf{929} & \textbf{929} & \textbf{929} \\
163 & 878 & \cite{MarshallPhd2011} & \textbf{934} & \textbf{934} & \textbf{934} & \textbf{934} \\
164 & 880 & \cite{MarshallPhd2011} & \textbf{940} & \textbf{940} & \textbf{940} & \textbf{940} \\
165 & 883 & \cite{MarshallPhd2011} & \textbf{946} & \textbf{946} & \textbf{946} & \textbf{946} \\
166 & 886 & \cite{MarshallPhd2011} & \textbf{953} & \textbf{953} & \textbf{953} & \textbf{953} \\
167 & 892 & \cite{AbajoEtAl2010new} & \textbf{958} & \textbf{958} & \textbf{958} & \textbf{958} \\
168 & 901 & \cite{AbajoEtAl2010new} & \textbf{964} & \textbf{965} & \textbf{965} & \textbf{965} \\
169 & 910 & \cite{AbajoEtAl2010new} & \textbf{970} & \textbf{971} & \textbf{971} & \textbf{971} \\
170 & 920 & \cite{AbajoEtAl2010new} & \textbf{977} & \textbf{978} & \textbf{978} & \textbf{978} \\
171 & 930 & \cite{AbajoEtAl2010new} & \textbf{984} & \textbf{984} & \textbf{984} & \textbf{984} \\
172 & 932 & \cite{AbajoEtAl2010new} & \textbf{990} & \textbf{992} & \textbf{992} & \textbf{992} \\
173 & 935 & \cite{AbajoEtAl2010new} & \textbf{997} & \textbf{998} & \textbf{998} & \textbf{998} \\
174 & 938 & \cite{AbajoEtAl2010new} & \textbf{1005} & \textbf{1005} & \textbf{1005} & \textbf{1005} \\
175 & 941 & \cite{AbajoEtAl2010new} & \textbf{1012} & \textbf{1012} & \textbf{1012} & \textbf{1012} \\
176 & 949 & \cite{MarshallPhd2011} & \textbf{1019} & \textbf{1020} & \textbf{1020} & \textbf{1020} \\
177 & 958 & \cite{MarshallPhd2011} & \textbf{1027} & \textbf{1027} & \textbf{1027} & \textbf{1027} \\
178 & 968 & \cite{MarshallPhd2011} & \textbf{1034} & \textbf{1035} & \textbf{1035} & \textbf{1035} \\
179 & 977 & \cite{MarshallPhd2011} & \textbf{1042} & \textbf{1042} & \textbf{1042} & \textbf{1042} \\
180 & 986 & \cite{MarshallPhd2011} & \textbf{1050} & \textbf{1050} & \textbf{1050} & \textbf{1050} \\
181 & 995 & \cite{MarshallPhd2011} & \textbf{1055} & \textbf{1056} & \textbf{1056} & \textbf{1056} \\
182 & 1004 & \cite{MarshallPhd2011} & \textbf{1061} & \textbf{1063} & \textbf{1063} & \textbf{1063} \\
183 & 1013 & \cite{MarshallPhd2011} & \textbf{1067} & \textbf{1069} & \textbf{1069} & \textbf{1069} \\
184 & 1022 & \cite{MarshallPhd2011} & \textbf{1073} & \textbf{1075} & \textbf{1076} & \textbf{1076} \\
185 & 1032 & \cite{MarshallPhd2011} & \textbf{1079} & \textbf{1082} & \textbf{1082} & \textbf{1082} \\
186 & 1042 & \cite{MarshallPhd2011} & \textbf{1087} & \textbf{1088} & \textbf{1088} & \textbf{1088} \\
187 & 1052 & \cite{MarshallPhd2011} & \textbf{1093} & \textbf{1094} & \textbf{1094} & \textbf{1094} \\
188 & 1062 & \cite{MarshallPhd2011} & \textbf{1100} & \textbf{1101} & \textbf{1101} & \textbf{1101} \\
189 & 1072 & \cite{MarshallPhd2011} & \textbf{1106} & \textbf{1107} & \textbf{1107} & \textbf{1107} \\
190 & 1082 & \cite{MarshallPhd2011} & \textbf{1112} & \textbf{1114} & \textbf{1114} & \textbf{1114} \\
191 & 1092 & \cite{MarshallPhd2011} & \textbf{1118} & \textbf{1120} & \textbf{1120} & \textbf{1120} \\
192 & 1102 & \cite{MarshallPhd2011} & \textbf{1125} & \textbf{1126} & \textbf{1126} & \textbf{1126} \\
193 & 1112 & \cite{MarshallPhd2011} & \textbf{1131} & \textbf{1132} & \textbf{1133} & \textbf{1133} \\
194 & 1122 & \cite{MarshallPhd2011} & \textbf{1138} & \textbf{1139} & \textbf{1139} & \textbf{1139} \\
195 & 1132 & \cite{MarshallPhd2011} & \textbf{1144} & \textbf{1145} & \textbf{1146} & \textbf{1146} \\
196 & 1142 & \cite{MarshallPhd2011} & \textbf{1152} & \textbf{1152} & \textbf{1152} & \textbf{1153} \\
\bottomrule
\end{tabular}
\caption{Lower bounds of $ex(n;\{C_3,C_4\})$ for $158 \leq n \leq 196$ from the literature and \cref{algo:compLowerBounds}. Bounds which improve upon the ones from the literature are marked in bold.}
\label{tab:results158to196}
\end{table}

\begin{table}[h]
\centering
\begin{tabular}{r r r |  r r r r }
\toprule
\multicolumn{1}{c}{$n$} & \multicolumn{6}{c}{$ex(n;\{C_3,C_4\})\geq$} \\
\cmidrule(lr){2-7}
 & & & \multicolumn{4}{c}{\cref{algo:compLowerBounds}} \\ 
 \cmidrule(lr){4-7}
 & Literature & Source & Up 1 & Down 1 & Up 2 & Down 2 \\ 
\midrule
197 & 1152 & \cite{MarshallPhd2011} & \textbf{1159} & \textbf{1159} & \textbf{1159} & \textbf{1160} \\
198 & 1163 & \cite{MarshallPhd2011} & \textbf{1166} & \textbf{1166} & \textbf{1166} & \textbf{1166} \\
199 & 1173 & \cite{MarshallPhd2011} & 1172 & 1173 & 1173 & 1173 \\
200 & 1184 & \cite{MarshallPhd2011} & 1179 & 1184 & 1184 & 1184 \\
201 & 1195 & \cite{MarshallPhd2011} & 1186 & 1195 & 1195 & 1195 \\
202 & 1206 & \cite{MarshallPhd2011} & 1193 & 1206 & 1206 & 1206 \\
203 & 1218 & \cite{MarshallPhd2011} & 1218 & - & 1218 & - \\
\bottomrule
\end{tabular}
\caption{Lower bounds of $ex(n;\{C_3,C_4\})$ for $197 \leq n \leq 203$ from the literature and \cref{algo:compLowerBounds}. Bounds which improve upon the ones from the literature are marked in bold.}
\label{tab:results197to203}
\end{table}

\end{document}